\documentclass[a4paper,reqno,12pt]{amsart}
\usepackage[utf8]{inputenc}
\usepackage[T1,T2A]{fontenc}
\usepackage[english]{babel}
\usepackage{amsmath}
\usepackage{amssymb}
\usepackage{amscd}
\usepackage{amsthm}
\usepackage{euscript}
\usepackage[all]{xy}
\newtheorem{proposition}{Proposition}[section]
\newtheorem{conjecture}[proposition]{Conjecture}
\newtheorem{problem}[proposition]{Problem}

\newtheorem{lemma}[proposition]{Lemma}
\newtheorem{theorem}[proposition]{Theorem}
\newtheorem{lettertheorem}{Theorem}

\theoremstyle{definition}
\newtheorem{definition}[proposition]{Definition}

\theoremstyle{remark}
\newtheorem{remark}[proposition]{Remark}

\DeclareMathOperator{\Char}{char}

\def\FF{{\mathbb F}}

\def\RRR{\mathfrak{R}}
\def\SSS{\mathfrak{S}}
\def\ggg{\mathfrak{g}}
\def\sss{\mathfrak{s}}

\def\zzz{\mathfrak{z}}
\def\gl{\mathfrak{gl}}
\def\sl{\mathfrak{sl}}

\renewcommand{\phi}{\varphi}
\renewcommand{\ge}{\geqslant}
\renewcommand{\le}{\leqslant}

\begin{document}
    \date{}
    \title[Uniqueness of addition in Lie rings $\gl_n(K)$ and $\sl_n(K)$]{Uniqueness of addition in Lie rings $\gl_n(K)$ and $\sl_n(K)$}
    \author{Gennadiy Sosnov}
    \address{Faculty of Computer Science, HSE University, Pokrovsky Boulevard 11, Moscow, 109028 Russia}
    \email{grsosnov@edu.hse.ru}
	\thanks{This article is an output of a research project (HSE-BR-2025-22) implemented as part of the Basic Research Program at HSE University.}
    
    \begin{abstract}
    We prove that for any Lie ring $\RRR$ and any commutator-preserving bijection $\alpha : \gl_n(K) \rightarrow \RRR$, the map $\alpha$ is additive on $\sl_n(K)$, where $n \ge 2$ and $K$ is an arbitrary field. Using this result we find criteria for Lie rings $\gl_n(K)$ and $\sl_n(K)$ to be unique addition. We also show that any commutator-preserving injection of Lie rings $\beta : \sl_2(K) \rightarrow \SSS$ is additive. This is the first result on additivity of commutator-preserving injections.
    \end{abstract}
    \maketitle


    \section{Introduction}
    \label{sec1}

    An associative ring $(R, +, \cdot)$ is called a unique addition ring, or UA-ring for short, if for any ring $(S, \oplus, *)$ every isomorphism $\alpha : R \rightarrow S$ of the multiplicative semigroups $(R, \cdot)$ and $(S, *)$ is an isomorphism of rings. This concept traces back to Rickart's works. In~\cite{Rickart1948} he proved that every Boolean ring is a UA-ring. Later, Stephenson in~\cite{Stephenson1969} showed that every matrix ring $M_n(R)$ is a UA-ring, where $n \ge 2$ and $R$ is an arbitrary associative ring with unit.

    However, this theory is not well suited for non-associative rings. The first results on uniqueness of addition for Lie rings were obtained by Arzhantsev. In~\cite{Arzhantsev2000} it was proven that the Lie algebra $\sl_2(K)$ is a UA-Lie algebra for any field $K$ with $\Char K \neq 2$. A year later, it was shown in~\cite{Arzhantsev2001} that every semisimple Lie algebra $\ggg$ over an algebraically closed field of characteristic zero is a UA-ring. Mayorova~\cite{Mayorova2019} proved that the uniqueness of addition holds in the case of Lie algebras of Chevalley type over rings with $1/2$ and $1/3$.
	
	Stronger results were obtained by Dolinar, Kuzma and Marovt~\cite{Dolinar}. They explicitly described the structure of commutator-preserving bijections from $M_n(K)$ to themselves. Namely, for any commutator-preserving bijection $\phi : M_n(K) \rightarrow M_n(K)$, where $K$ is a field, there exists an invertible matrix $T \in M_n(K)$, a field automorphism $\sigma : K \rightarrow K$ and a function $\psi : M_n(K) \rightarrow K$, where $\psi(A) = 0$ for all matrices of trace zero such that
	\begin{enumerate}
		\item for $n \ge 3$ and $K$ with at least $2^{n - 1}$ elements, either
		\begin{equation*}
			\phi(A) = TA^{\sigma}T^{-1} + \psi(A)E \quad \text{for all} \quad A \in M_n(K),
		\end{equation*}
		or
		\begin{equation*}
			\phi(A) = -T(A^{\sigma})^{\top}T^{-1} + \psi(A)E \quad \text{for all} \quad A \in M_n(K);
		\end{equation*}
		\item for $n = 2$ and $\Char K \neq 2$,
		\begin{equation*}
			\phi(A) = TA^{\sigma}T^{-1} + \psi(A)E \quad \text{for all} \quad A \in M_n(K).
		\end{equation*}
	\end{enumerate}
	Here $E$ is the identity matrix and $A^{\sigma}$ is the matrix obtained by applying $\sigma$ to all matrix entries.

	This article is devoted to the case of matrix Lie algebras $\sl_n(K)$ and $\gl_n(K)$ over arbitrary fields. Let $\FF_q$ denote a finite field with $q$ elements.

	\begin{lettertheorem}\label{theorem:when_glIsUA}
		Let $n \ge 2$. Then the ring $\gl_n(K)$ is a UA-ring if and only if $n$ is odd and $K = \FF_2$.
	\end{lettertheorem}
	
	This result answers the question posed in~\cite[Problem 1]{arzhantsev2025uniqueness}.
	
	\begin{lettertheorem}\label{theorem:when_slIsUA}
		Let $n \ge 2$. Then the ring $\sl_n(K)$ is not a UA-ring if and only if $\Char K = 2$ and $n = 2$.
	\end{lettertheorem}

	We also consider a natural generalization of the concept of UA-rings to the case of injective maps.

	\begin{lettertheorem}\label{theorem:injectionsForsl2}
		Let $K$ be a field with $\Char K \neq 2$ and $\SSS$ be a Lie ring. Then any commutator-preserving injection $\beta : \sl_2(K) \rightarrow \SSS$ is additive.
	\end{lettertheorem}

	We prove Theorem~A, Theorem~B and Theorem~C in Sections~\ref{sec3}, \ref{sec5}, \ref{sec4}, respectively. Section~\ref{sec2} contains basic definitions and preliminary results.

	The author is grateful to Ivan Arzhantsev for posing the problem and useful suggestions throughout this work.

	\section{Preliminary results}
	\label{sec2}

	We start with the main definition for this work.

	\begin{definition}
		A Lie ring $\RRR$ is called a \textit{unique addition Lie ring} or a \textit{UA-Lie ring} if for any Lie ring $\SSS$ and any commutator-preserving bijection $\alpha : \RRR \rightarrow \SSS$, the map $\alpha$ is additive, i.e. the condition $\alpha([a, b]) = [\alpha(a), \alpha(b)]$ for all $a, b \in \RRR$ implies $\alpha(a + b) = \alpha(a) + \alpha(b)$.
	\end{definition}

	Let us define $C(a) = \lbrace b \in \RRR \mid [a, b] = 0 \rbrace$ to be the centralizer of an element $a$ and $Z(\RRR) = \lbrace x \in \RRR \mid \forall y \in \RRR : [x, y] = 0 \rbrace$ to be the center of a Lie ring $\RRR$. Also, let us fix a commutator-preserving bijection $\alpha : \RRR \rightarrow \SSS$. Clearly, $\alpha(0) = \alpha([0, 0]) = [\alpha(0), \alpha(0)] = 0$. Now we consider the properties of $\alpha$ in terms of the center and centralizers.

	\begin{lemma}\label{lemma:bijectionSavesCenter}
		The equality $\alpha(Z(\RRR)) = Z(\SSS)$ holds.
	\end{lemma}

	\begin{proof}
		Let $z \in Z(\RRR)$. Then for all $x \in \RRR$ we have $0 = \alpha([x, z]) = [\alpha(x), \alpha(z)]$, therefore $\alpha(z) \in Z(\SSS)$. Hence, $\alpha(Z(\RRR)) \subseteq Z(\SSS)$. Applying this argument to the commutator-preserving bijection $\alpha^{-1}$, we obtain the required equality.
	\end{proof}

	Let us recall a definition that simplifies the study of UA-Lie rings.

	\begin{definition}
		Let $\alpha : \RRR \rightarrow \SSS$ be any commutator-preserving map of Lie rings. The map $\Delta : \RRR \times \RRR \rightarrow \SSS$ that sends each pair $(a, b)$ to the element $\alpha(a + b) - \alpha(a) - \alpha(b)$ is called the \textit{defect} of $\alpha$.
	\end{definition}

	In~\cite[Lemma~1]{Arzhantsev2001} it was proved that for any $a, b \in \RRR$ we have
	\begin{equation*}
		\Delta(a, b) \in \alpha(C(a) \cap C(b)).
	\end{equation*}
	Moreover, a stronger statement holds.

	\begin{lemma}\label{lemma:intersectionOfCenters}
		For any $a, b \in \RRR$ we have
		\begin{equation*}
			\Delta(a, b) \in \alpha(Z(C(a)) \cap Z(C(b))).
		\end{equation*}
	\end{lemma}

	\begin{proof}
		Let $\alpha(a + b) = \alpha(a) + \alpha(b) + \alpha(c)$. It is enough to show that $c \in Z(C(a)) \cap Z(C(b))$. Let $x$ be any element in $C(a)$. Since
		\begin{equation*}
			\alpha([b, x]) = \alpha([a + b, x]) = [\alpha(a), \alpha(x)] + [\alpha(b), \alpha(x)] + [\alpha(c), \alpha(x)] = \alpha([b, x]) + \alpha([c, x]),
		\end{equation*}
		we have $\alpha([c, x]) = 0$, therefore $[c, x] = 0$ and $c \in Z(C(a))$. Similarly, we obtain ${c \in Z(C(b))}$.
	\end{proof}

	A natural generalization of a UA-Lie ring is the concept of an $iUA$-ring.

	\begin{definition}
		A Lie ring $\RRR$ is called an \textit{iUA-Lie ring} if for any Lie ring $\SSS$ and any commutator-preserving injection $\beta : \RRR \rightarrow \SSS$, the map $\beta$ is additive.
	\end{definition}

	It is easy to show that Lemmas~\ref{lemma:bijectionSavesCenter} and \ref{lemma:intersectionOfCenters} fail in the case of injections. Nevertheless, the defect still has some interesting properties. Let us fix a commutator-preserving injection of Lie rings $\beta : \RRR \rightarrow \SSS$ and its defect $\Delta$.
	
	\begin{lemma}\label{lemma:JacobiForAdditivity}
		Let $x, y \in \RRR$ be elements such that there exist $a, b, c \in \RRR$ with 
		\begin{equation*}
			x = [a, [b, c]] \quad \text{and} \quad y = [b, [c, a]].
		\end{equation*}
		Then $\Delta(x, y) = 0$.
	\end{lemma}

	\begin{proof}
		By the Jacobi identity we have
		\begin{equation*}
			x + y = [a, [b, c]] + [b, [c, a]] = [[a, b], c],
		\end{equation*}
		so
		\begin{equation*}
			\beta(x + y) = [[\beta(a), \beta(b)], \beta(c)] = [\beta(a), [\beta(b), \beta(c)]] + [\beta(b), [\beta(c), \beta(a)]] = \beta(x) + \beta(y).
		\end{equation*}
		Therefore, $\Delta(x, y) = 0$.
	\end{proof}

	\begin{lemma}\label{lemma:commutatorForAdditivity}
		For any $x, y, z \in \RRR$ we have
		\begin{equation*}
			[\beta(x), \Delta(y, z)] = \Delta([x, y], [x, z]).
		\end{equation*}
		In particular, $\Delta([x, y], [x, z]) = 0$ provided $\Delta(y, z) = 0$.
	\end{lemma}

	\begin{proof}
		By using the commutator-preserving property, we obtain
		\begin{equation*}
			[\beta(x), \Delta(y, z)] = [\beta(x), \beta(y + z)] - [\beta(x), \beta(y)] - [\beta(x), \beta(z)] = \Delta([x, y], [x, z]).
		\end{equation*}
	\end{proof}

	\begin{lemma}\label{lemma:sumOfDefects}
		For any $x, y, z \in \RRR$ we have
		\begin{equation*}
			\Delta(x + y, z) + \Delta(x, y) = \Delta(x, y + z) + \Delta(y, z).
		\end{equation*}
	\end{lemma}

	\begin{proof}
		The left hand side is
		\begin{equation*}
			\beta(x + y + z) - \beta(x + y) - \beta(z) + \beta(x + y) - \beta(x) - \beta(y) = \beta(x + y + z) - \beta(x) - \beta(y) - \beta(z),
		\end{equation*}
		and the right hand side is
		\begin{equation*}
			\beta(x + y + z) - \beta(x) - \beta(y + z) + \beta(y + z) - \beta(y) - \beta(z) = \beta(x + y + z) - \beta(x) - \beta(y) - \beta(z).
		\end{equation*}
	\end{proof}

	\section{UA-properties for the general linear Lie algebra}
	\label{sec3}

	Let $E$ be the identity matrix and $E_{ij}$ be a matrix unit. It is shown in~\cite[Corollary~2, Remark~2]{arzhantsev2025uniqueness} that the Lie algebra $\gl_n(K)$ is not a UA-Lie ring if either $K \neq \FF_2$ or $K = \FF_2$ and $n$ is even. For the convenience of the reader, we provide a short proof of this result.

	\begin{proposition}\label{proposition:whenglIsNotUA}
		The algebra $\gl_n(K)$ is not a UA-Lie ring if either $K \neq \FF_2$ or $K = \FF_2$ and $n$ is even.
	\end{proposition}

	\begin{proof}
		If the field $K$ consists of more than $2$ elements, then we can construct the bijection $\alpha : \gl_n(K) \rightarrow \gl_n(K)$.
		\begin{equation*}
			\alpha(X + tE) = X + \psi(t)E, \quad \text{for all} \quad X \in \sl_n(K), t \in K,
		\end{equation*}
		where $\psi$ is an arbitrary non-additive field bijection. It is obvious that the map $\alpha$ preserves the commutator and is not additive.

		On the other hand, if $K = \FF_2$ and $n$ is even, take a matrix $A$ with non-zero trace. Clearly, $A \not\in [\gl_n(K), \gl_n(K)]$ and $A + E \not\in [\gl_n(K), \gl_n(K)]$, therefore the map $\alpha : \gl_n(K) \rightarrow \gl_n(K)$ that swaps the elements $A$ and $A + E$ preserves the commutator but is not additive.
	\end{proof}

	Let us fix a commutator-preserving bijection $\alpha : \gl_n(K) \rightarrow \SSS$, where $n \ge 3$ (the case $n = 2$ is considered in Proposition~\ref{proposition:injectionsForgl2}) with defect $\Delta$. Although $\alpha$ is not additive in general, it is almost additive. More precisely, the following proposition holds.

	\begin{proposition}\label{proposition:glIsAlmostUA}
		Let $K$ be a field and $\alpha : \gl_n(K) \rightarrow \SSS$ be a commutator-preserving bijection. Then $\alpha$ is additive on $\sl_n(K)$, i.e. $\alpha(A + B) = \alpha(A) + \alpha(B)$ for any $A, B \in \sl_n(K)$.
	\end{proposition}

	Before proving this proposition, we need the following lemmas. It is important to remember that $n \ge 3$; some of these proofs do not work in the case $n = 2$.

	\begin{lemma}
		For any matrix unit $E_{ij} \in \gl_n(K)$ we have $Z(C(E_{ij})) = \langle E_{ij}, E \rangle$.
	\end{lemma}

	\begin{proof}
		If $i = j$ then $C(E_{ii}) = \langle E_{ii} \rangle \oplus \gl_{n - 1}(K)$. Therefore, $Z(C(E_{ii})) = \langle E_{ii}, E \rangle$. Otherwise, if $i \neq j$ then for all $s \neq j$ and $t \neq i$ we have $E_{st} \in C(E_{ij})$. Since $n > 2$, these conditions ensure that $Z(C(E_{ij})) = \langle E_{ij}, E \rangle$. 
	\end{proof}
	
	\begin{lemma}\label{lemma:ZCForDiag}
		Let $A \in \gl_n(K)$ be a diagonal matrix. Then the set $Z(C(A))$ consists of diagonal matrices.
	\end{lemma}

	\begin{proof}
		Let $X \in Z(C(A))$. Clearly, $E_{ii} - E_{jj} \in C(A)$ for any $i \neq j$. Therefore, the equality $[X, E_{ii} - E_{jj}] = 0$ holds, so we have $X_{ik} = X_{ki} = X_{jk} = X_{kj} = 0$ for all $k \neq i$ and $k \neq j$. Since $n > 2$, the previous observation shows that $X$ is a diagonal matrix.
	\end{proof}

	It is known that the center $Z(\gl_n(K))$ is the set of all scalar matrices. Let us denote it as $\zzz$.

	\begin{lemma}\label{lemma:additivityForBothDiag}
		Let $A, B \in \gl_n(K)$ be two diagonal matrices. Then $\Delta(A, B) \in \alpha(\zzz)$.
	\end{lemma}

	\begin{proof}
		Since the matrix units $E_{ii}$ span the space of diagonal matrices, it is sufficient to show that $\Delta(A, \lambda E_{ii}) \in \alpha(\zzz)$ for all $i$ and $\lambda \in K$ by Lemma~\ref{lemma:sumOfDefects}. If $Z(C(A)) \cap Z(C(\lambda E_{ii})) = \zzz$, then our lemma follows from Lemma~\ref{lemma:intersectionOfCenters}. Otherwise, $Z(C(\lambda E_{ii})) \subseteq Z(C(A))$, so $[A, \lambda E_{ii}] = 0$. Let $C = E_{ik}$, where $k \neq i$. Then $[\lambda E_{ii}, A + C] = [\lambda E_{ii}, C] \neq 0$. Therefore, the equality $Z(C(\lambda E_{ii})) \cap Z(C(A + C)) = \zzz$ holds. Since $A + \lambda E_{ii}$ and $A$ are diagonal matrices, the sets $Z(C(A + \lambda E_{ii}))$ and $Z(C(A))$ consist of diagonal matrices. Hence,
		\begin{equation*}
			Z(C(A + \lambda E_{ii})) \cap Z(C(C)) = Z(C(A)) \cap Z(C(C)) = \zzz.
		\end{equation*}
		Applying Lemmas~\ref{lemma:intersectionOfCenters} and \ref{lemma:sumOfDefects}, we obtain
		\begin{equation*}
			\Delta(A, \lambda E_{ii}) = \Delta(\lambda E_{ii}, A + C) + \Delta(A, C) - \Delta(A + \lambda E_{ii}, C) \in \alpha(\zzz).
		\end{equation*}
	\end{proof}

	\begin{lemma}\label{lemma:almostAdditivityForDiag}
		Let $A, B \in \gl_n(K)$ and $B$ be a diagonal matrix. Then $\Delta(A, B) \in \alpha(\zzz)$.
	\end{lemma}

	\begin{proof}
		Again, it is sufficient to prove the statement in the case $B = \lambda E_{ii}$. If $A$ is a diagonal matrix, then the lemma follows from Lemma~\ref{lemma:additivityForBothDiag}. Therefore, without loss of generality, we can assume that $A$ is not a diagonal matrix. Let $s, t$ be such different indices that the inequality $A_{st} \neq 0$ holds. If $Z(C(A)) \cap Z(C(\lambda E_{ii})) = \zzz$ then the lemma follows from Lemma~\ref{lemma:intersectionOfCenters}. Otherwise, $Z(C(\lambda E_{ii})) \subseteq Z(C(A))$, so $[A, \lambda E_{ii}] = 0$. Since $t \neq i$, we have $[\lambda E_{ii}, A + E_{ti}] \neq 0$, $[A, E_{ti}] \neq 0$ and $[A + \lambda E_{ii}, E_{ti}] \neq 0$. Therefore,
		\begin{equation*}
			Z(C(A + \lambda E_{ii})) \cap Z(C(E_{ti})) = Z(C(\lambda E_{ii})) \cap Z(C(A + E_{ti})) = Z(C(A)) \cap Z(C(E_{ti})) = \zzz.
		\end{equation*}
		Hence,
		\begin{equation*}
			\Delta(A, \lambda E_{ii}) = \Delta(\lambda E_{ii}, A + E_{ti}) + \Delta(A, E_{ti}) - \Delta(A + \lambda E_{ii}, E_{ti}) \in \alpha(\zzz).
		\end{equation*}
	\end{proof}

	\begin{lemma}\label{lemma:almostAdditive}
		Let $A, B \in \gl_n(K)$, then $\Delta(A, B) \in \alpha(\zzz)$.
	\end{lemma}

	\begin{proof}
		It is sufficient to prove the lemma for $B = \lambda E_{ij}$. The case $i = j$ was proved in the previous lemma. So, without loss of generality, we can assume that $i \neq j$. Let $C$ be a diagonal matrix such that $[A + C, B] \neq 0$. Then
		\begin{equation*}
			\Delta(A + B, C), \Delta(B, A + C), \Delta(A, C) \in \alpha(\zzz).
		\end{equation*}
		Applying Lemma~\ref{lemma:sumOfDefects}, we obtain
		\begin{equation*}
			\Delta(A, B) = \Delta(B, A + C) + \Delta(A, C) - \Delta(A + B, C) \in \alpha(\zzz).
		\end{equation*}
	\end{proof}

	\begin{lemma}\label{lemma:additivityForDiag}
		The defect $\Delta(A, B) = 0$ for any diagonal matrices $A$ and $B$ in $\sl_n(K)$.
	\end{lemma}

	\begin{proof}
		Let $J$ be a nilpotent Jordan block, so $J = \sum_{i = 1}^{n - 1} E_{i, i + 1}$. The image of the adjoint operator of $J$ contains the space of all diagonal matrices of trace zero. Therefore, for any diagonal matrices $A, B$ there exist matrices $C, D$ such that $A = [J, C]$ and $B = [J, D]$. From Lemma~\ref{lemma:almostAdditive} it follows that $\Delta(C, D) \in \alpha(\zzz)$, hence
		\begin{equation*}
			\Delta(A, B) = [\alpha(J), \Delta(C, D)] = 0.
		\end{equation*}
	\end{proof}

	Let us fix a standard basis $\mathcal{B} = \lbrace E_{ij} \mid i \neq j \rbrace \cup \lbrace E_{11} - E_{ii} \mid i \neq 1 \rbrace$ of $\sl_n(K)$. The key observation of the whole proof is the theorem proved by Albert and Muckenhoupt~\cite{albert_muckenhoupt1957}.

	\begin{theorem}\label{theorem:commutatorWidthForgl}
		Any matrix in $\sl_n(K)$ is a commutator of two matrices in $\gl_n(K)$.
	\end{theorem}

	\begin{lemma}\label{lemma:additivityForProportional}
		Let $X \in \sl_n(K)$ and $\lambda_1, \lambda_2 \in K$. Then $\Delta(\lambda_1 X, \lambda_2 X) = 0$.
	\end{lemma}

	\begin{proof}
		Since $X = [A, B]$ for some matrices $A, B \in \gl_n(K)$ and $\Delta(\lambda_1 B, \lambda_2 B) \in \alpha(\zzz)$, it follows from Lemma~\ref{lemma:commutatorForAdditivity} that
		\begin{equation*}
			\Delta(\lambda_1 X, \lambda_2 X) = \Delta([A, \lambda_1 B], [A, \lambda_2 B]) = [\alpha(A), \Delta(\lambda_1 B, \lambda_2 B)] = 0.
		\end{equation*}
	\end{proof}

	\begin{lemma}\label{lemma:basisDecomposition}
		Let $X$ be any matrix in $\sl_n(K)$ and $X = \sum_{i, j} \beta_{ij} E_{ij} + \sum_{i} \beta_{i} (E_{11} - E_{ii})$ be the decomposition of $X$ in the basis $\mathcal{B}$. Then
		\begin{equation*}
			\alpha(X) = \sum_{i, j} \alpha(\beta_{ij} E_{ij}) + \sum_{i} \alpha(\beta_{i} (E_{11} - E_{ii})).
		\end{equation*}
	\end{lemma}

	\begin{proof}
		Let $X = [A, B]$ for some matrices $A, B \in \gl_n(K)$ and
		\begin{equation*}
			A = \sum_{i, j} a_{ij} E_{ij}, \quad B = \sum_{i, j} b_{ij} E_{ij}
		\end{equation*}
		be their decompositions. Then
		\begin{equation*}
			\alpha(X) = [\alpha(A), \alpha(B)] = \sum_{i,j,k,s} \alpha([a_{ij}E_{ij}, b_{ks}E_{ks}]).
		\end{equation*}
		Since $[a_{ij}E_{ij}, b_{ks}E_{ks}]$ is proportional to a basis element in $\mathcal{B}$ or is equal to a diagonal matrix, it is sufficient to prove additivity for diagonal matrices, which we already have done in Lemma~\ref{lemma:additivityForDiag}.
 	\end{proof}

	\begin{proof}[Proof of Proposition~\ref{proposition:glIsAlmostUA}]
		Let us consider two arbitrary matrices $A, B \in \sl_n(K)$ with decompositions in $\mathcal{B}$
		\begin{equation*}
			A = \sum_{i, j} \beta_{ij} E_{ij} + \sum_{i} \beta_{i} (E_{11} - E_{ii})
		\end{equation*}
		and
		\begin{equation*}
			B = \sum_{i, j} \gamma_{ij} E_{ij} + \sum_{i} \gamma_{i} (E_{11} - E_{ii}).
		\end{equation*}
		Then the sum $A + B$ has the following decomposition
		\begin{equation*}
			A + B = \sum_{i, j} (\beta_{ij} + \gamma_{ij}) E_{ij} + \sum_{i} (\beta_{i} + \gamma_{i}) (E_{11} - E_{ii}),
		\end{equation*}
		so applying Lemma~\ref{lemma:basisDecomposition} we obtain
		\begin{multline*}
			\alpha(A + B) = \sum_{i, j} \alpha((\beta_{ij} + \gamma_{ij}) E_{ij}) + \sum_{i} \alpha((\beta_{i} + \gamma_{i}) (E_{11} - E_{ii})) = \\ = \sum_{i, j} \left[ \alpha(\beta_{ij} E_{ij}) + \alpha(\gamma_{ij} E_{ij}) \right] + \sum_{i} \left[ \alpha(\beta_{i}(E_{11} - E_{ii})) + \alpha(\gamma_i (E_{11} - E_{ii})) \right] = \alpha(A) + \alpha(B).
		\end{multline*}
	\end{proof}

	With this proposition we can obtain criteria for Lie ring $\gl_n(K)$ to be unique addition.

	\begin{proof}[Proof of Theorem~\ref{theorem:when_glIsUA}]
		The negative results were provided in Proposition~\ref{proposition:whenglIsNotUA}, so it is sufficient to prove that the ring $\gl_{2n + 1}(\FF_2)$ is a UA-Lie ring for any $n \ge 1$. Let $\alpha$ be a commutator-preserving bijection from $\gl_{2n + 1}(\FF_2)$ to the Lie ring $\SSS$. Since $\alpha(E) + \alpha(E) \neq \alpha(E)$ and $|Z(\SSS)| = 2$, we have $\alpha(E) + \alpha(E) = 0 = \alpha(E + E)$. Therefore, we only need to prove that $\alpha(A + E) = \alpha(A) + \alpha(E)$ for any matrix $A$ in $\sl_{2n + 1}(\FF_2)$. Assume the contrary, then $\alpha(A + E) = \alpha(A) + \alpha(E) + \alpha(T)$ for some matrix $T \in Z(C(A)) \cap Z(C(E)) = \zzz$ and $T \neq 0$. Thus, $T = E$ and $\alpha(A + E) = \alpha(A)$. We obtain a contradiction, so the map $\alpha$ is additive.
	\end{proof}

	\section{Commutator-preserving injections}
	\label{sec4}

	Let $\beta : \RRR \rightarrow \SSS$ be a commutator-preserving injection between two Lie rings with defect~$\Delta$. Our goal is to determine when $\beta$ is additive.

	The article~\cite{arzhantsev2025uniqueness} states the following conjecture.

	\begin{conjecture}\label{conjecture:injectiveMaps}
		Every semisimple Lie algebra is an iUA-Lie ring.
	\end{conjecture}

	The same article proves the following negative result~\cite[Lemma~9]{arzhantsev2025uniqueness}.

	\begin{lemma}
		Let $K$ be a field with at least $3$ elements and $\ggg$ be a Lie algebra over $K$ such that $\ggg \neq [\ggg, \ggg]$. Then there exists a Lie algebra $\sss$ and a non-additive commutator-preserving injection $\beta : \ggg \rightarrow \sss$.
	\end{lemma}

	This lemma can be easily generalized.

	\begin{lemma}
		Let $\RRR$ be a Lie ring with at least $2$ elements and $\RRR \neq \lbrace [x, y] \mid x, y \in \RRR \rbrace$. Then the ring $\RRR$ is not an iUA-Lie ring.
	\end{lemma}

	\begin{proof}
		Let $\SSS = \RRR \oplus \FF_3$, where $\FF_3$ is an abelian Lie ring with exactly three elements. Also, let $y$ be a fixed element in $\RRR \setminus \lbrace [x, y] \mid x, y \in \RRR \rbrace$. Then the desired injection $\beta : \RRR \rightarrow \SSS$ is defined as
		\begin{equation*}
			\beta(x) = 
			\begin{cases}
				(x, 0), \quad x \neq y \\
				(x, 1), \quad x = y
			\end{cases}.
		\end{equation*}
		The injection $\beta$ preserves the commutator because
		\begin{equation*}
			\beta([a, b]) = ([a, b], 0) = [\beta(a), \beta(b)], \quad \text{for all} \quad a, b \in \RRR.
		\end{equation*}
		However,
		\begin{equation*}
			(y + y, 0) = \beta(y + y) \neq \beta(y) + \beta(y) = (y + y, 2),
		\end{equation*}
		hence $\beta$ is not additive.
	\end{proof}

	Therefore, every iUA-Lie ring has bracket width equal to $1$. For semisimple Lie algebras the question of maximal bracket width is open (see~\cite{kunyavskii2021bracket}). However, the article~\cite{kunyavskii2021bracket} constructs an infinite-dimensional simple Lie algebra with bracket width strictly greater than $1$, so there are simple Lie algebras that are not iUA-Lie rings.

	We will prove Conjecture~\ref{conjecture:injectiveMaps} in the case of $\ggg = \sl_2(K)$, where $\Char K \neq 2$. For the proof, we need to fix some notation. Let $e, f, h$ be the standard basis of $\sl_2(K)$, so

	\begin{equation*}
		e = 
		\begin{pmatrix}
			0 & 1 \\
			0 & 0
		\end{pmatrix}, \quad
		h =
		\begin{pmatrix}
			1 & 0 \\
			0 & -1
		\end{pmatrix}, \quad
		f = 
		\begin{pmatrix}
			0 & 0 \\
			1 & 0
		\end{pmatrix}.
	\end{equation*}

	\begin{proposition}\label{proposition:injectionsForgl2}
		Let $K$ be a field with arbitrary characteristic. Then any commutator-preserving injection of Lie rings $\beta : \gl_2(K) \rightarrow \SSS$ is additive on $\sl_2(K)$.
	\end{proposition}

	To do this, we need to prove several lemmas.

	\begin{lemma}\label{lemma:additivityForEandH}
		For any $\lambda_1, \lambda_2 \in K$ we have $\Delta(\lambda_1 h, \lambda_2 e) = 0$ and $\Delta(\lambda_1 h, \lambda_2 f) = 0$.
	\end{lemma}

	\begin{proof}
		The case $\lambda_1 = 0$ or $\lambda_2 = 0$ is obvious. Therefore, we will consider the case $\lambda_1, \lambda_2 \neq 0$. Applying Lemma~\ref{lemma:JacobiForAdditivity} to the matrices
		\begin{equation*}
			A = 
			\begin{pmatrix}
				0 & \lambda_2 \\
				0 & 0
			\end{pmatrix}, \quad
			B = 
			\begin{pmatrix}
				1 & 0 \\
				0 & 0
			\end{pmatrix}, \quad
			C = 
			\begin{pmatrix}
				1 & 0 \\
				-\frac{\lambda_1}{\lambda_2} & 0
			\end{pmatrix},
		\end{equation*}
		we obtain $\Delta(\lambda_1 h, \lambda_2 e) = \Delta([A, [B, C]], [B, [C, A]]) = 0$. Since the identity
		\begin{equation*}
			[A^{\top}, [B^{\top}, C^{\top}]] = [A, [B, C]]^{\top}
		\end{equation*}
		holds, we get the equality $\Delta(\lambda_1 h, \lambda_2 f) = 0$.
	\end{proof}

	\begin{lemma}\label{lemma:additivityForE}
		For any $\lambda_1, \lambda_2 \in K$ we have
		\begin{equation*}
			\Delta(\lambda_1 e, \lambda_2 e) = 0 \quad \text{and} \quad \Delta(\lambda_1 f, \lambda_2 f) = 0.
		\end{equation*}
	\end{lemma}
	
	\begin{proof}
		Let
		\begin{equation*}
			A = 
			\begin{pmatrix}
				1 & 0 \\
				0 & 0
			\end{pmatrix}, \quad
			B = 
			\begin{pmatrix}
				1 & 1 \\
				0 & 0
			\end{pmatrix}, \quad
			C = 
			\begin{pmatrix}
				-\lambda_1 - \lambda_2 & -\lambda_2 \\
				0 & 0
			\end{pmatrix}.
		\end{equation*}
		Then $[A, [B, C]] = \lambda_1 e$ and $[B, [C, A]] = \lambda_2 e$. Applying Lemma~\ref{lemma:JacobiForAdditivity}, we obtain the equality $\Delta(\lambda_1 e, \lambda_2 e) = 0$. The equality $\Delta(\lambda_1 f, \lambda_2 f) = 0$ is obtained by transposition.
	\end{proof}

	\begin{lemma}\label{lemma:additivityForHandMinusH}
		For any $\lambda \in K$ we have $\Delta(\lambda h, -\lambda h) = 0$.
	\end{lemma}

	\begin{proof}
		Let
		\begin{equation*}
			A = 
			\begin{pmatrix}
				0 & 1 \\
				0 & 0
			\end{pmatrix}, \quad
			B = 
			\begin{pmatrix}
				0 & 0 \\
				\lambda & 0
			\end{pmatrix}, \quad
			C = 
			\begin{pmatrix}
				1 & 0 \\
				0 & 0
			\end{pmatrix}.
		\end{equation*}
		Then $[A, [B, C]] = \lambda h$ and $[B, [C, A]] = -\lambda h$. Applying Lemma~\ref{lemma:JacobiForAdditivity}, we obtain the required result.
	\end{proof}

	\begin{remark}\label{remark:phiIsNotSurjective}
		In fact, for any field $K$ with $\Char K \neq 2$ and any $\lambda_1, \lambda_2 \in K \setminus \lbrace 0 \rbrace$ with $\lambda_1 + \lambda_2 \neq 0$, there are no matrices $A, B, C \in \gl_2(K)$ such that
		\begin{equation*}
			\lambda_1 h = [A, [B, C]], \quad \lambda_2 h = [B, [C, A]].
		\end{equation*}
		Because otherwise $\lambda_3 h = -(\lambda_1 + \lambda_2) h = [C, [A, B]]$ and since $\gl_2(K) = \sl_2(K) \oplus \langle E \rangle$ we can assume that $A, B, C \in \sl_2(K)$. If $\varkappa$ is the Killing form on $\sl_2(K)$, then
		\begin{align*}
			\varkappa(\lambda_1 h, A) &= \varkappa([A, [B, C]], A) = \varkappa([C, B], [A, A]) = 0 \Rightarrow A \in \langle h \rangle^{\perp}, \\
			\varkappa(\lambda_2 h, B) &= \varkappa([B, [C, A]], B) = \varkappa([A, C], [B, B]) = 0 \Rightarrow B \in \langle h \rangle^{\perp}, \\
			\varkappa(\lambda_3 h, C) &= \varkappa([C, [A, B]], C) = \varkappa([B, A], [C, C]) = 0 \Rightarrow C \in \langle h \rangle^{\perp}.
		\end{align*}
		But $\langle h \rangle^{\perp} = \langle e, f \rangle$, therefore $[B, C] \in \langle h \rangle$ and $[A, [B, C]] \in \langle e, f \rangle$. We obtain a contradiction.
	\end{remark}

	\begin{lemma}\label{lemma:additivityForEandF}
		For any $\lambda_1, \lambda_2 \in K$ we have $\Delta(\lambda_1 e, \lambda_2 f) = 0$.
	\end{lemma}

	\begin{proof}
		Let
		\begin{equation*}
			A = 
			\begin{pmatrix}
				1 & 0 \\
				\lambda_2 & 0
			\end{pmatrix}, \quad
			B = 
			\begin{pmatrix}
				1 & -\lambda_1 \\
				0 & 0
			\end{pmatrix}, \quad
			C = 
			\begin{pmatrix}
				1 & 0 \\
				0 & 0
			\end{pmatrix}.
		\end{equation*}
		Then $[A, [B, C]] = -\lambda_1 \lambda_2 h + \lambda_1 e$ and $[B, [C, A]] = \lambda_1\lambda_2 h + \lambda_2 f$. Therefore,
		\begin{equation*}
			0 = \Delta(-\lambda_1 \lambda_2 h + \lambda_1 e, \lambda_1\lambda_2 h + \lambda_2 f) = \beta(\lambda_1 e + \lambda_2 f) - \beta(-\lambda_1 \lambda_2 h + \lambda_1 e) - \beta(\lambda_1\lambda_2 h + \lambda_2 f).
		\end{equation*}
		Applying Lemmas~\ref{lemma:additivityForEandH} and \ref{lemma:additivityForHandMinusH}, we obtain
		\begin{multline*}
			\beta(\lambda_1 e + \lambda_2 f) - \beta(-\lambda_1\lambda_2 h + \lambda_1 e) - \beta(\lambda_1\lambda_2 h + \lambda_2 f) = \\ = \beta(\lambda_1 e + \lambda_2 f) - \beta(-\lambda_1\lambda_2 h) - \beta(\lambda_1 e) - \beta(\lambda_1\lambda_2 h) - \beta(\lambda_2 f) = \Delta(\lambda_1 e, \lambda_2 f).
		\end{multline*}
		Therefore, $\Delta(\lambda_1 e, \lambda_2 f) = 0$.
	\end{proof}

	\begin{lemma}\label{lemma:additivityForH}
		For any $\lambda_1, \lambda_2 \in K$ we have $\Delta(\lambda_1 h, \lambda_2 h) = 0$.
	\end{lemma}

	\begin{proof}
		Let
		\begin{equation*}
			A = 
			\begin{pmatrix}
				0 & \lambda_2 \\
				-\lambda_1 & 0
			\end{pmatrix}.
		\end{equation*}
		Applying Lemmas~\ref{lemma:additivityForEandF} and \ref{lemma:commutatorForAdditivity}, we obtain
		\begin{equation*}
			\Delta(\lambda_1 h, \lambda_2 h) = \Delta([A, e], [A, f]) = [\beta(A), \Delta(e, f)] = 0.
		\end{equation*}
	\end{proof}

	\begin{lemma}\label{lemma:additivityForTriangular}
		For any $\lambda_1, \lambda_2, \mu_1, \mu_2 \in K$ we have
		\begin{equation*}
			\Delta(\lambda_1 h + \lambda_2 e, \mu_1 h + \mu_2 f) = 0.
		\end{equation*}
	\end{lemma}

	\begin{proof}
		If $\lambda_2 = 0$, then applying Lemmas~\ref{lemma:additivityForEandH} and \ref{lemma:additivityForH}, we obtain the required result. The case $\mu_2 = 0$ is considered similarly. So we can assume $\lambda_2 \neq 0$ and $\mu_2 \neq 0$. Let
		\begin{equation*}
			A = 
			\begin{pmatrix}
				1 & -\frac{\mu_1}{\mu_2} \\
				-\frac{\lambda_1}{\lambda_2} & 0
			\end{pmatrix}.
		\end{equation*}
		Then applying Lemmas~\ref{lemma:additivityForEandF} and \ref{lemma:commutatorForAdditivity}, we obtain
		\begin{equation*}
			\Delta(\lambda_1 h + \lambda_2 e, \mu_1 h + \mu_2 f) = \Delta([A, \lambda_2 e], [A, -\mu_2 f]) = [\beta(A), \Delta(\lambda_2 e, -\mu_2 f)] = 0.
		\end{equation*}
	\end{proof}

	\begin{proof}[Proof of Proposition~\ref{proposition:injectionsForgl2}]
		Let
		\begin{equation*}
			A = 
			\begin{pmatrix}
				a_1 & b_1 \\
				c_1 & -a_1
			\end{pmatrix}, \quad
			B =
			\begin{pmatrix}
				a_2 & b_2 \\
				c_2 & -a_2
			\end{pmatrix}.
		\end{equation*}
		It follows from Lemma~\ref{lemma:additivityForTriangular} that
		\begin{equation*}
			\beta(A + B) = \beta(a_1 h + (b_1 + b_2)e) + \beta(a_2 h + (c_1 + c_2)f).
		\end{equation*}
		Therefore, applying Lemma~\ref{lemma:additivityForEandH}, we obtain
		\begin{equation*}
			\beta(a_1 h + (b_1 + b_2)e) + \beta(a_2 h + (c_1 + c_2)f) = \beta(a_1 h + b_1 e) + \beta(b_2 e) + \beta(a_2 h + c_2 f) + \beta(c_1 f).
		\end{equation*}
		Again, it follows from Lemma~\ref{lemma:additivityForTriangular} that 
		\begin{equation*}
			\beta(a_1 h + b_1 e) + \beta(c_1 f) = \beta(A), \quad \beta(a_2 h + c_2 f) + \beta(b_2 e) = \beta(B),
		\end{equation*}
		so $\beta(A + B) = \beta(A) + \beta(B)$.
	\end{proof}

	With this proposition we can easily prove Theorem~\ref{theorem:injectionsForsl2}.

	\begin{proof}[Proof of Theorem~\ref{theorem:injectionsForsl2}]
		Let $\beta : \sl_2(K) \rightarrow \SSS$ be a commutator-preserving injection of Lie rings. Let $\SSS' = \SSS \oplus Kz$ be the Lie ring with bracket $[\SSS, Kz] = 0$ and $\beta' : \gl_2(K) \rightarrow \SSS'$ be the map sending $A + kE$ to $\beta(A) + kz$, where $A \in \sl_2(K)$ and $k \in K$. The map $\beta'$ preserves the commutator, therefore $\beta'$ is additive on $\sl_2(K)$. Since $\beta'|_{\sl_2(K)} = \beta$, we obtain the required result.
	\end{proof}

	Using a computer, one can check that the map
	\begin{equation*}
		\phi : \RRR \times \RRR \times \RRR \rightarrow [\RRR, \RRR] \times [\RRR, \RRR], \quad (a, b, c) \mapsto ([a, [b, c]], [b, [c, a]])
	\end{equation*}
	is surjective for $\RRR = \gl_3(\FF_2)$. Therefore, applying Lemma~\ref{lemma:JacobiForAdditivity} and modifying the proof of Theorem~\ref{theorem:injectionsForsl2}, we can show that the ring $\sl_3(\FF_2)$ is an iUA-Lie ring. From Remark~\ref{remark:phiIsNotSurjective} it follows that the map $\phi$ is not surjective for $\RRR = \gl_2(K)$, where $\Char K \neq 2$.

	\begin{problem}
		Is $\phi$ surjective for $\RRR = \gl_n(K)$, where $n \ge 3$?
	\end{problem}

	\section{UA-properties for the special linear Lie algebra}
	\label{sec5}

	First, let us prove the negative result.

	\begin{lemma}\label{lemma:whenSlIsNotUA}
		The Lie ring $\sl_2(K)$ is not a UA-Lie ring if $\Char K = 2$.
	\end{lemma}

	\begin{proof}
		In this case $[\sl_2(K), \sl_2(K)] = \langle E \rangle$, where $E$ is the identity matrix. Therefore, we can construct a bijection $\alpha : \sl_2(K) \rightarrow \sl_2(K)$ in the following way
		\begin{equation*}
			\alpha(x) = 
			\begin{cases}
				e + E, & x = e \\
				e, & x = e + E \\
				x, & \text{otherwise}
			\end{cases}.
		\end{equation*}
		Clearly, it preserves the commutator and is not additive.
	\end{proof}

	For $\Char K \neq 2$ we know from Theorem~\ref{theorem:injectionsForsl2} that the algebra $\sl_2(K)$ is an iUA-Lie ring, so it is definitely a UA-Lie ring. So we only need to prove that the Lie algebra $\sl_n(K)$ is a UA-Lie ring for every $n \ge 3$ and an arbitrary field $K$. Let the map $\alpha : \sl_n(K) \rightarrow \SSS$ be a commutator-preserving bijection of Lie rings and $\Delta$ be the defect of $\alpha$.

	\begin{lemma}\label{lemma:slIsUAwhenCharNDevideN}
		Let $\Char K \nmid n$ and $n \ge 3$. Then the Lie algebra $\sl_n(K)$ is a UA-Lie ring.
	\end{lemma}

	\begin{proof}
		Since $\gl_n(K) = \sl_n(K) \oplus \langle E \rangle$, we can extend the map $\alpha$ to $\alpha' : \gl_n(K) \rightarrow \SSS'$ by defining $\SSS'$ as $\SSS \oplus Kz$, where $z$ is a central element and
		\begin{equation*}
			\alpha'(A + kE) = \alpha(A) + kz.
		\end{equation*}
		Clearly, the map $\alpha'$ preserves the commutator. Applying Proposition~\ref{proposition:glIsAlmostUA}, we obtain that $\alpha'$ is additive on $\sl_n(K)$. Since $\alpha'|_{\sl_n(K)} = \alpha$, we get the required result.
	\end{proof}

	So we can assume that $\Char K \mid n$. It is known from the work of Thompson~\cite{Thompson1966} that any matrix $A \in \sl_n(K)$ can be written as $A = [X, Y]$ for some $X, Y \in \sl_n(K)$, except when $\Char K = 2$ and $n = 2$. So the bracket width of the Lie algebra $\sl_n(K)$ is equal to $1$ for an arbitrary field $K$ and $n \ge 3$. Let us denote $Z(\sl_n(K))$ as $\zzz = \langle E \rangle$. Now it is clear from the proof of Proposition~\ref{proposition:glIsAlmostUA} that we only need to prove that $\Delta(A, B) \in \zzz$ for any $A, B \in \sl_n(K)$.
	
	\begin{lemma}
		Let $E_{ij} \in \sl_n(K)$ be a matrix unit, where $i \neq j$. Then $Z(C(E_{ij})) = \langle E_{ij}, E \rangle$.
	\end{lemma}

	\begin{proof}
		For all $s \neq j$ and $t \neq i$ we have $E_{st} \in C(E_{ij})$. Since $n > 2$, these conditions ensure that $Z(C(E_{ij})) = \langle E_{ij}, E \rangle$. 
	\end{proof}
	
	\begin{lemma}\label{lemma:ZCForDiagInSl}
		Let $A \in \sl_n(K)$ be a diagonal matrix. Then the set $Z(C(A))$ consists only of diagonal matrices.
	\end{lemma}

	\begin{proof}
		Let $X \in Z(C(A))$. Clearly, $E_{ii} - E_{jj} \in C(A)$ for any $i \neq j$. Therefore, the equality $[X, E_{ii} - E_{jj}] = 0$ holds, so we have $X_{ik} = X_{ki} = X_{jk} = X_{kj} = 0$ for all $k \neq i$ and $k \neq j$. Since $n > 2$, the previous observation shows us that $X$ is a diagonal matrix.
	\end{proof}

	\begin{lemma}
		For any $i \neq j$ the following holds
		\begin{equation*}
			Z(C(E_{ii} - E_{jj})) = \langle E_{ii} - E_{jj}, E \rangle.
		\end{equation*}
	\end{lemma}

	\begin{proof}
		From Lemma~\ref{lemma:ZCForDiagInSl} we know that if $X \in Z(C(E_{ii} - E_{jj}))$, then $X$ is diagonal. Let $s$ and $t$ be two indices not equal to $i$ and $j$. Since $[X, E_{st}] = 0$, we conclude that $X_{ss} = X_{tt}$ for all such $s$ and $t$. Assume that $\dim Z(C(E_{ii} - E_{jj})) \ge 3$. Then we can find linearly independent matrices
		\begin{equation*}
			E_{ii} - E_{jj}, E, Y \in Z(C(E_{ii} - E_{jj})).
		\end{equation*}
		Therefore, we can assume that $Y_{ss} = 0$ for all $s \neq i$. Since $\operatorname{tr}(Y) = 0$, we conclude that $Y = 0$. We obtain a contradiction, so $\dim Z(C(E_{ii} - E_{jj})) \le 2$. Since the inclusion $\langle E_{ii} - E_{jj}, E \rangle \subseteq Z(C(E_{ii} - E_{jj}))$ holds, we get the required result.
	\end{proof}

	\begin{lemma}\label{lemma:additivityForBothDiagInSl}
		Let $A, B \in \sl_n(K)$ be two diagonal matrices. Then $\Delta(A, B) \in \alpha(\zzz)$.
	\end{lemma}

	\begin{proof}
		Since the matrices $E_{ii} - E_{jj}$ span the space of diagonal matrices, it is sufficient to prove the lemma for $B = \lambda (E_{ii} - E_{jj})$, where $i \neq j$ and $\lambda \in K$. If $Z(C(A)) \cap Z(C(B)) = \zzz$, then the lemma follows from Lemma~\ref{lemma:intersectionOfCenters}. Otherwise, $Z(C(B)) \subseteq Z(C(A))$, hence the equality $[A, B] = 0$ holds. Let $C = E_{ik}$, where $k \neq i$ and $k \neq j$. Then $[B, A + C] \neq 0$, therefore $Z(C(B)) \cap Z(C(A + C)) = \zzz$. Since $A + B$ and $A$ are diagonal matrices, we conclude that the sets $Z(C(A + B))$ and $Z(C(A))$ consist of diagonal matrices. Hence,
		\begin{equation*}
			Z(C(A + B)) \cap Z(C(C)) = Z(C(A)) \cap Z(C(C)) = \zzz.
		\end{equation*}
		Applying Lemmas~\ref{lemma:intersectionOfCenters} and \ref{lemma:sumOfDefects}, we obtain
		\begin{equation*}
			\Delta(A, B) = \Delta(B, A + C) + \Delta(A, C) - \Delta(A + B, C) \in \alpha(\zzz).
		\end{equation*}
	\end{proof}

	\begin{lemma}\label{lemma:almostAdditivityForDiagInSl}
		Let $A, B \in \sl_n(K)$ and $B$ be a diagonal matrix. Then $\Delta(A, B) \in \alpha(\zzz)$.
	\end{lemma}

	\begin{proof}
		Again, it is sufficient to prove the statement in the case $B = \lambda (E_{ii} - E_{jj})$. If $A$ is a diagonal matrix, then the lemma follows from Lemma~\ref{lemma:additivityForBothDiag}. Therefore, without loss of generality, we can assume that $A$ is not a diagonal matrix. Let $A_{st} \neq 0$, where $s \neq t$. If $Z(C(A)) \cap Z(C(B)) = \zzz$, then the lemma follows from Lemma~\ref{lemma:intersectionOfCenters}. Otherwise, $Z(C(B)) \subseteq Z(C(A))$, so $[A, B] = 0$.

		\textit{Case 1.} Let $s \in \lbrace i, j \rbrace$, $t \in \lbrace i, j \rbrace$ and $k \not\in \lbrace i, j \rbrace$. Since $[B, A + E_{ks}] = [B, E_{ks}] \neq 0$, $[A, E_{ks}] \neq 0$ and $[A + B, E_{ks}] \neq 0$, we have
		\begin{equation*}
			Z(C(A + B)) \cap Z(C(E_{ks})) = Z(C(B)) \cap Z(C(A + E_{ks})) = Z(C(A)) \cap Z(C(E_{ks})) = \zzz.
		\end{equation*}
		Hence,
		\begin{equation*}
			\Delta(A, B) = \Delta(B, A + E_{ks}) + \Delta(A, E_{ks}) - \Delta(A + B, E_{ks}) \in \alpha(\zzz).
		\end{equation*}

		\textit{Case 2.} Let $s \neq i$ and $s \neq j$. Since $[B, A + E_{is}] = [B, E_{is}] \neq 0$, $[A, E_{is}] \neq 0$ and $[A + B, E_{is}] \neq 0$, we have
		\begin{equation*}
			Z(C(A + B)) \cap Z(C(E_{is})) = Z(C(B)) \cap Z(C(A + E_{is})) = Z(C(A)) \cap Z(C(E_{is})) = \zzz.
		\end{equation*}
		Hence,
		\begin{equation*}
			\Delta(A, B) = \Delta(B, A + E_{is}) + \Delta(A, E_{is}) - \Delta(A + B, E_{is}) \in \alpha(\zzz).
		\end{equation*}
		
		\textit{Case 3.} Let $t \neq i$ and $t \neq j$. Since $[B, A + E_{ti}] = [B, E_{ti}] \neq 0$, $[A, E_{ti}] \neq 0$ and $[A + B, E_{ti}] \neq 0$, we have
		\begin{equation*}
			Z(C(A + B)) \cap Z(C(E_{ti})) = Z(C(B)) \cap Z(C(A + E_{ti})) = Z(C(A)) \cap Z(C(E_{ti})) = \zzz.
		\end{equation*}
		Hence,
		\begin{equation*}
			\Delta(A, B) = \Delta(B, A + E_{ti}) + \Delta(A, E_{ti}) - \Delta(A + B, E_{ti}) \in \alpha(\zzz).
		\end{equation*}
	\end{proof}

	Let us again fix a standard basis $\mathcal{B} = \lbrace E_{ij} \mid i \neq j \rbrace \cup \lbrace E_{11} - E_{ii} \mid i \neq 1 \rbrace$ of $\sl_n(K)$.

	\begin{lemma}\label{lemma:almostAdditiveInSl}
		The defect $\Delta(A, B) \in \alpha(\zzz)$ for any $A, B \in \sl_n(K)$.
	\end{lemma}

	\begin{proof}
		It is sufficient to prove this lemma for matrix $B$ that is proportional to an element in $\mathcal{B}$. The case $B = \lambda(E_{ii} - E_{jj})$ was proved in the previous lemma. So, without loss of generality, we can assume that $B = \lambda E_{ij}$, where $i \neq j$. Let $C$ be a diagonal matrix such that $[A + C, B] \neq 0$. Then
		\begin{equation*}
			\Delta(A + B, C), \Delta(B, A + C), \Delta(A, C) \in \alpha(\zzz).
		\end{equation*}
		Applying Lemma~\ref{lemma:sumOfDefects}, we obtain
		\begin{equation*}
			\Delta(A, B) = \Delta(B, A + C) + \Delta(A, C) - \Delta(A + B, C) \in \alpha(\zzz).
		\end{equation*}
	\end{proof}

	Since the algebra $\sl_n(K)$ has commutator width $1$, we can repeat the proof of Proposition~\ref{proposition:glIsAlmostUA}, but for the convenience of the reader we summarize the whole proof of Theorem~\ref{theorem:when_slIsUA}.

	\begin{proof}[Proof of Theorem~\ref{theorem:when_slIsUA}]
		From Lemma~\ref{lemma:whenSlIsNotUA} and Theorem~\ref{theorem:injectionsForsl2} we obtain that the Lie algebra $\sl_2(K)$ is a UA-Lie ring if and only if $\Char K \neq 2$. From Lemma~\ref{lemma:slIsUAwhenCharNDevideN} we conclude that $\sl_n(K)$ is a UA-Lie ring if $\Char K \nmid n$.
		
		Now let $n \ge 3$, $\Char K \mid n$ and $\alpha : \sl_n(K) \rightarrow \SSS$ be a commutator-preserving bijection of Lie rings with defect $\Delta$. From Lemma~\ref{lemma:almostAdditiveInSl} we obtain that $\Delta(A, B) \in Z(\sl_n(K))$ for any $A, B \in \sl_n(K)$. Similarly to what we did in the proof of Lemma~\ref{lemma:additivityForDiag}, we obtain that $\Delta(A, B) = 0$ for any diagonal matrices $A, B \in \sl_n(K)$. Since the Lie ring $\sl_n(K)$ has bracket width $1$, we get that $\Delta(\lambda_1 X, \lambda_2 X) = 0$ for any $\lambda_1, \lambda_2 \in K$ and $X \in \sl_n(K)$. Therefore, $\alpha$ preserves decompositions with respect to the basis $\mathcal{B}$, so $\alpha$ is additive.
	\end{proof}


\end{document}